\documentclass[12pt, reqno]{amsart}
\usepackage{amsmath,amssymb,amsthm,graphicx,mathrsfs,url} % Required for inserting images
\usepackage{amsmath,amssymb,amsthm,tikz}
\usetikzlibrary{decorations.pathreplacing}
\usepackage[toc,page]{appendix}

\usepackage[margin=1in]{geometry}
\usepackage[hidelinks]{hyperref}
\usepackage{comment}

\usepackage[normalem]{ulem}

\newcommand{\R}{\mathbb R}

\newtheorem{theorem}{Theorem}[section]
\newtheorem{lemma}[theorem]{Lemma}
\newtheorem{proposition}[theorem]{Proposition}

\newtheorem{corollary}[theorem]{Corollary}

\theoremstyle{remark}
\newtheorem{remark}[theorem]{Remark}

\newcommand{\eps}{\varepsilon}
\newcommand{\la}{\lambda}

\newcommand{\tr}{\mathrm{tr}}

\numberwithin{equation}{section}

\title[Monotonicity of topological entropy along the Ricci flow]{Monotonicity of topological entropy along the Ricci flow near a hyperbolic metric}
\author{Karen Butt, Alena Erchenko, Tristan Humbert}
\date{}

\begin{document}
\begin{abstract}
In 2004, Manning showed that the topological entropy of the geodesic flow of a closed surface of non-constant negative curvature is strictly decreasing along the normalized Ricci flow, and he asked if an analogous result holds in higher dimensions for metrics in a neighborhood of a hyperbolic metric. 
In this paper, we affirmatively answer this question. Namely, we show that the topological entropy of the geodesic flow of a closed Riemannian manifold that carries a hyperbolic metric is indeed strictly decreasing along the normalized Ricci flow starting from a metric of variable negative sectional curvature sufficiently close to the hyperbolic metric.
\end{abstract}
\maketitle

\section{Introduction}
%\subsection{Statement of main result}

The Ricci flow is a natural deformation on the space of metrics on a closed manifold that was introduced by Hamilton \cite{Hamilton} as a tool to uniformize a metric.
In this paper, we study how the \emph{topological entropy} of the geodesic flow, a numerical invariant which quantifies the total complexity of a dynamical system's orbit structure (see Section~\ref{section: topological entropy}), changes as we evolve the underlying metric by the \emph{normalized Ricci flow} (see \eqref{eq:NRF-intro}).

Suppose $M$ is a closed manifold of dimension at least 3 which admits a \emph{hyperbolic metric}, i.e., a metric of constant negative sectional curvature (unique up to isometry and rescaling by Mostow rigidity \cite{mostow}).
Then, by work of Besson--Courtois--Gallot \cite{BCGGAFA}, the topological entropy of any non-hyperbolic metric on $M$ (of the same total volume) is strictly larger than that of the hyperbolic metric. By work of Ye \cite{Ye}, the normalized Ricci flow starting from a metric sufficiently $C^2$-close to a hyperbolic metric is defined for all positive times and converges exponentially to a hyperbolic metric. 

Thus, a natural question, first raised by Manning \cite[Question 5]{Man_Ricci}, is whether the topological entropy of the starting metric strictly decreases to that of the hyperbolic metric in this setting. 
Our main theorem provides an affirmative answer to this question. 
\begin{theorem}
    \label{MainTheo} Let $(M^n,g_{\mathrm{hyp}})$ be a closed hyperbolic manifold of dimension $n\geq 3$. There exist $N=N(n)\in \mathbb N$ and $\epsilon>0$ such that for any non-hyperbolic metric $g$ satisfying $\|g~-~g_{\mathrm{hyp}}\|_{C^N}<\epsilon$, the topological entropy is strictly decreasing along the normalized Ricci flow $(g_t)_{t\geq 0}$ starting from $g_0 = g$.
\end{theorem}

\begin{remark}
 In our setting, the topological entropy coincides with the \emph{volume entropy}, which measures the exponential growth rate (in $R$) of the volume of a geodesic ball of radius $R$ in the universal cover; see \cite{manningannals}.
\end{remark}

We recall that 
the \emph{normalized Ricci flow} (NRF) $(g_t)_{t\geq 0}$ starting from $g$ is defined by the differential equation
\begin{equation}
\label{eq:NRF-intro}
\frac{\partial}{\partial{t}}g_t=-2\mathrm{Ric}_t +\frac{2}{n}\frac{\int_M\mathrm{Scal}_t d\mathrm{vol}_t }{\mathrm{Vol}_g(M)}g_t, \quad g_0=g,
\end{equation} 
where $\mathrm{Ric}_t$ is the Ricci tensor, $\mathrm{Scal}_t$ is the scalar curvature, $d\mathrm{vol}_t$ is the volume form associated to $g_t$, and $\mathrm{Vol}_g(M)$ is the total volume of $M$ with respect to $g$. (See Section \ref{sec:curv} for the precise definitions.) The average scalar curvature term on the right-hand side of \eqref{eq:NRF-intro} ensures the total volume of $M$ with respect to $g_t$ remains constant and is the reason for the use of the word ``normalized".

It is important to note that the behavior of the Ricci flow in dimensions $n \geq 3$ is quite different from the case $n = 2$. In dimension 2, the normalized Ricci flow can be used to prove that every metric on a closed surface is conformally equivalent to a metric of constant curvature, a classical fact which follows from the uniformization theorem for Riemann surfaces.
Indeed, the NRF starting from any metric $g$ on a closed surface is defined for all times and converges to a metric of constant curvature; moreover, the NRF in dimension $2$ preserves negative Gaussian curvature and the conformal class of the starting metric. 
However, in dimensions $n \geq 3$, the NRF does \emph{not} preserve the conformal class of the starting metric and does \emph{not} preserve negative sectional curvature. What is more, the NRF starting from a negatively curved metric in higher dimensions need not converge to an Einstein metric.
In fact, Farrell and Ontaneda \cite{FO} showed that given any $\epsilon>0$, there exists a hyperbolic manifold $(M,g_{\rm hyp})$ and a negatively curved $(1-\epsilon)$-pinched metric $g$ on $M$ such that the NRF starting from $g$ does \emph{not} converge in the $C^2$-topology to a negatively curved Einstein metric. 

\subsection{Related results}
The analogue of Theorem \ref{MainTheo} in dimension 2 follows from work of Manning \cite[Theorem 1]{Man_Ricci}. 
In fact, Manning proved the stronger result that $h_{\rm top}(g)$ is strictly decreasing along the NRF starting from \emph{any} metric $g$ of variable negative curvature, not necessarily close to $g_{\rm hyp}$. 
On the other hand, if the starting metric $g$ is not negatively curved, the topological entropy can increase along the NRF.
For instance, Jane \cite{Ja} constructed examples of smooth
metrics on the 2-sphere and the 2-torus which have zero topological entropy and for which the entropy strictly increases along the NRF. 
Also, Vasii \cite{vasii} showed that the volume growth entropy, which coincides with the topological entropy for metrics of nonpositive sectional curvature \cite{manningannals}, does not decrease along the (unnormalized) Ricci flow $(g_t)_{t \geq 0}$, provided $g_t$ does not collapse.

Manning asked if a similar monotonicity result to \cite[Theorem 1]{Man_Ricci} holds in higher dimensions. More precisely, he asked if $t \mapsto h_{\rm top}(g_t)$ is strictly decreasing for $g_0$ in a neighborhood of a hyperbolic metric $g_{\rm hyp}$ \cite[Question 5]{Man_Ricci} (see also \cite{Thom}). The caveat that $g_0$ is close to $g_{\rm hyp}$ is natural in light of both the aforementioned positive results of Ye \cite{Ye} on the convergence of the Ricci flow starting from such a $g_0$ as well as the results of Farrell and Ontaneda \cite{FO} on the non-convergence of the NRF for metrics with ``almost constant" negative curvature.

A partial answer to \cite[Question 5]{Man_Ricci} was previously obtained by Thompson \cite{Thom}.
He showed that the derivative of the topological entropy along the NRF at $t=0$ is negative when starting from a metric with negative sectional curvature which is conformally equivalent (but not isometric) to $g_{\mathrm{hyp}}$, in any dimension. In particular, this implies that the topological entropy decreases along the NRF for small times $t>0$. To compare, the monotonicity result in Theorem \ref{MainTheo} holds for any metric near $g_{\mathrm{hyp}}$ (not necessarily in the same conformal class) and holds for any time $t\geq 0.$ We note, however, that Thompson does not suppose that the starting metric $g$ is close to the hyperbolic metric $g_{\rm hyp}$, whereas we do. We remark that while Thompson and Manning's proofs coincide when $n = 2$, our method is distinct (see Section~\ref{strat} below for more details).

In dimension $2$, the NRF coincides with another important geometric flow, the \emph{Yamabe flow}. 
In light of this, a higher dimensional generalization of \cite[Theorem 1]{Man_Ricci} can also be interpreted as studying the variation of $h_{\rm top}(g)$ under the Yamabe flow \cite[Question 5]{Man_Ricci}.
Related questions were considered by 
Suarez-Serrato and Tapie \cite{ST}. They defined two flows, the \emph{increasing/decreasing curvature-normalized Yamabe flow} in dimension $n\geq 3$ and studied their convergences. Moreover, they established entropy bounds (with corresponding rigidity statements) along their flows. In particular, they obtained a monotonicity result for the volume growth entropy along these flows \cite[Proof of Theorem 24]{ST}.
 
A related question of Manning \cite[Question 3]{Man_Ricci} asks about the monotonicity of other invariants along the NRF. In a previous paper, the authors together with Mitsutani \cite{BEHM} showed that the \emph{Liouville entropy} is increasing along the NRF starting from a $\frac{1}{6}$-pinched negatively curved metric on a surface and the \emph{mean root curvature} is also increasing along the NRF starting from any negatively curved metric. It is not clear whether these monotonicity results  extend to higher dimensions, even in a neighborhood of a hyperbolic metric. Indeed, contrary to the topological entropy, which is globally minimized at hyperbolic metrics, the Liouville entropy need not even have local extrema at such metrics when $n \geq 3$. Hyperbolic metrics are always critical points for the Liouville entropy, but they can be saddle points \cite{Fla,Maubon}.

\subsection{Stability estimate}

As in \cite{Man_Ricci,Thom}, we compute the derivative of the topological entropy $h_{\rm{top}}$ along the normalized Ricci flow at $t=0$ and show that it is negative in a neighborhood of $g_{\rm hyp}$. Recall that for a metric $g$ near $g_{\rm{hyp}}$, there is a unique smooth diffeomorphism homotopic to the
identity $\phi_g$ such that $\phi_g^*g-g_{\rm{hyp}}$ is \emph{divergence-free}, see Lemma~\ref{slice lemma}. We show the following estimate on the derivative of $h_{\rm top}$ along the NRF.
\begin{theorem}
\label{theo:2}Let $(M^n,g_{\mathrm{hyp}})$ be a closed hyperbolic manifold of dimension $n\geq 3$. There exist $N(n)\in \mathbb N$, $\epsilon>0$ and $C>0$ such that for any metric $g$ satisfying $\|g-g_{\mathrm{hyp}}\|_{C^N}<\epsilon$, if we let $(g_t)_{t\geq 0}$ be the normalized Ricci flow starting from $g$, then there is $c_g>0$ such that
\begin{equation}
    \label{eq:stability}
    \frac{d}{dt}h_{\rm{top}}(g_t)|_{t=0}\leq -C\|\phi_g^*g-c_gg_{\rm hyp}\|_{H^{-1/2}}^2,
\end{equation}
where $H^{-1/2}$ is the $L^2$-Sobolev space of order $-1/2$. In particular, if $\|g-g_{\rm hyp}\|_{C^N}<\epsilon$ and $g$ is not isometric to some rescalling of $g_{\rm hyp}$, then $\tfrac{d}{dt}h_{\rm{top}}(g_t)|_{t=0}<0.$ 
\end{theorem}

\begin{remark}
    Inspection of the proof of Theorem \ref{theo:2} shows that one can take $N(n)=\tfrac 32 n+17.$ 
\end{remark}

\begin{remark}
    Using an interpolation argument, one can replace the Sobolev norm $\Vert \cdot \Vert^2_{H^{-1/2}}$ with $\Vert \cdot \Vert_{C^0}^{\nu}$ for some explicit $\nu>0$ depending on $N$, and hence on $n$.
\end{remark}
We note that \cite[Theorem A]{BahuaudGuentherIsenberg} shows that there is some $C^{N+2}$-neighborhood $\mathcal U$ of $g_{\rm hyp}$ such that the (non-normalized) Ricci flow starting from any $g\in \mathcal U$ stays $C^N$-close to the Ricci flow starting from $g_{\rm hyp}$ (also called a Ricci soliton). Since the normalized Ricci flow can be obtained from the Ricci flow by an explicit time and space reparameterization, see, for instance, \cite[Chapter 9.1]{CK}, we obtain the existence of a $C^{N+2}$-neighborhood $\mathcal U'$ of $g_{\rm hyp}$ such that any NRF starting in $\mathcal U'$ stays $C^N$-close to $g_{\rm hyp}.$

In particular, Theorem \ref{theo:2}  implies Theorem \ref{MainTheo} (with possibly a slightly larger $N$ in Theorem \ref{MainTheo}). Moreover, we obtain a \emph{stability estimate} \eqref{eq:stability} which further gives a quantitative estimate on the derivative. We note that Manning and Thompson's arguments \cite{Man_Ricci,Thom} do not provide such estimates on the derivative.

We obtain estimate \eqref{eq:stability} using techniques from \emph{microlocal analysis}; see Section \ref{strat} for an outline of the argument. The use of microlocal tools to solve (local) rigidity problems in hyperbolic dynamics has become more common in recent years \cite{GL19,GuKnLef,GLP,Hum,BEHLW}.

\subsection{Strategy of the proof}\label{strat}
As mentioned above, our overall approach to proving Theorem~\ref{theo:2} is the same as in \cite{Man_Ricci, Thom}. Namely, we use a formula of Katok--Knieper--Weiss
\cite{KKW} for the first derivative of the topological entropy \eqref{eq:KKW} and plug in the equation for the normalized Ricci flow \eqref{eq:NRF-intro} to obtain
\begin{equation*}\label{intro:deriv}
    \left.\frac{d}{dt}\right|_{t=0}h_{\rm{top}}(g_t) =  -\frac{h_{\mathrm{top}}(g_0)}{2}\left(\int_{S^{0}M} -2{\rm {Ric}}_{0} (v) \, d \mu_{g_0} + \frac{2}{n {\rm {Vol}}(M)} \int_{M}\mathrm{Scal}_{0}d\mathrm{vol}_0\right),
\end{equation*}
where $\mu_{g_0}$ is the \emph{measure of maximal entropy} (also known as \emph{Bowen--Margulis measure}) of the geodesic flow of $g_0$ (see Section~\ref{section: topological entropy} for the definition).
%Plugging in the equation for the normalized Ricci flow \eqref{eq:NRF-intro} yields an expression of $\tfrac{d}{dt}|_{t=0}h_{\rm{top}}(g_t)$ as a functional of $g=g_0$. 
The strategy consists in studying the functional 
$g\mapsto \tfrac{d}{dt}|_{t=0}h_{\rm{top}}(g_t)$ with $g_0=g$
 in a neighborhood of $g_{\rm hyp}$ and showing it is negative.

Manning and Thompson's proofs rely on comparing the total scalar curvature term\linebreak $\mathcal S(g_t):=\int_M \mathrm{Scal}_t d\mathrm{vol}_t$ to  $\mathcal S(g_{\rm hyp})$.
In dimension 2, Manning uses that $\mathcal S(g_t)$ is constant  as a consequence of the Gauss--Bonnet theorem. In higher dimensions, Thompson's argument relies on the inequality $\mathcal S(g_t)\geq \mathcal S(g_{\rm hyp})$, which holds if $g_t$ is conformal to $g_{\rm hyp}$. However, from \cite[Proposition 4.55]{Bes}, we obtain that $\mathcal S(g)<\mathcal S(g_{\rm hyp})$ for some nearby trace-free and divergence-free perturbation of $g_{\rm hyp}$ when $n>2.$
As pointed out in \cite{Thom}, a new strategy is thus required to show negativity of $\tfrac{d}{dt}|_{t=0}h_{\rm{top}}(g_t)$ in greater generality when $n>2.$

Instead of using a zeroth-order approach to show the functional $g\mapsto \tfrac{d}{dt}|_{t=0}h_{\rm{top}}(g_t)$ is negative, our strategy is to use a second-order approach, that is, we compute the Hessian of this functional at $g_{\rm hyp}$ (or more precisely a closely related functional we call $\Phi$, see \eqref{eq:Phi}). 
While Manning and Thompson's arguments use the (global) inequality $h_{\rm top}(g_t) \geq h_{\rm top}(g_{\rm hyp})$ \cite{Ka, BCGGAFA}, our argument instead relies on Flaminio's computation of the Hessian of $h_{\rm top}$ at $g_{\rm hyp}$ (which implies that $h_{\rm top}$ has a local minimum at $g_{\rm hyp}$) \cite[Corollary 1.3.5]{Fla}; see also \cite{Pol}.

In Section \ref{sec:1stproperties}, we first check that $\Phi$ is a smooth functional in $g$ and that $g_{\rm hyp}$ is a zero and a critical point of $\Phi$. 
The core of the argument is to compute the Hessian of $\Phi$ at $g_{\rm hyp}$ and this is done in Proposition \ref{propHess1}. Since the first derivative of the topological entropy involves the Bowen--Margulis measure, our computation requires us to differentiate this measure. We solve this subtle problem using the computation of the Hessian of the topological entropy at real hyperbolic metrics due to Flaminio \cite{Fla}; see Proposition \ref{lemma:Flaminio}. Once the Hessian is simplified, we use Bochner/Weitzenböck inequalities to get a lower bound on the Hessian evaluated along divergence-free directions by a \emph{variance term}; see Proposition \ref{prop:var}.

We note that positivity of the Hessian is not enough to show that $\Phi$ is positive in a neighborhood of $g_0$. Hence, the last step of the proof is to deduce the stability estimate \eqref{eq:stability} from the lower bound given in Proposition \ref{prop:var}, which we do in Section \ref{sec:micro}. To this end, we use results from {microlocal analysis}, more precisely, the properties of the \emph{generalized X-ray transform} introduced by Guillarmou \cite{Gu} and used by Guillarmou-Lefeuvre \cite{GL19} (see also the related work of Guillarmou-Knieper-Lefeuvre \cite{GuKnLef}) to solve the local rigidity of the marked length spectrum. Note that the third-named author used a similar scheme to prove Katok's entropy conjecture near real and complex hyperbolic metrics \cite{Hum}. 

\subsection*{Acknowledgements}
We thank Colin Guillarmou and Thibault Lefeuvre for comments on an earlier draft of the paper.
Part of this paper was written while the authors attended an American Institute
of Mathematics (AIM) SQuaRE workshop in September 2025. We would like to thank AIM
for this opportunity and for their hospitality during our stay. 
K.B. was supported by NSF grant
DMS-2402173.
A.E. was supported by NSF grant DMS-2247230 and DMS-2552860. T.H was supported by the European Research Council (ERC) under
the European Union’s Horizon 2020 research and innovation programme (Grant agreement no. 101162990 — ADG).

\section{Preliminaries}

Let $(M,g)$ be a closed manifold equipped with a metric $g$ of negative sectional curvature. 
Let $d\mathrm{vol}_{g}$ denote the volume form associated to $g$ and $\mathrm{Vol}_{g}(M)=\int_M d\mathrm{vol}_g$ its total volume.
Let $S^{g}M:=\{(x,v)\in TM\,|\, \|v\|_{g}=1\}$ be its unit tangent bundle. We denote by $(\phi_g^t)_{t\in \mathbb R}$ the geodesic flow generated by $g$ on $S^gM.$ The generator of the flow $(\phi_g^t)_{t\in \mathbb R}$ is denoted by $X_g.$
\subsection{Topological entropy}\label{section: topological entropy}
In this paper, we study the topological entropy $h_{\rm top}(g)$ of the geodesic flow $\phi_g^t$ of $(M, g)$, equivalently, the entropy of the time-one map $\phi^1_g$.
The topological entropy of a dynamical system is a numerical invariant that quantifies the total exponential complexity of the system's orbit structure. 
More precisely, it measures the exponential rate of growth (in $t$) of the number of distinguishable orbit segments of length $t$. See, for instance, \cite[Chapter 4]{FishHas} for a more detailed account. For $g$ negatively curved (and, more generally, Anosov), it follows from work of Margulis \cite{margulis} and Bowen \cite{Bow72} that
\begin{equation}\label{eq:top}
h_{\rm top}(g) = \lim_{T \to \infty} \frac{\log \# \{ \gamma \, | \, l_{g}(\gamma) \leq T \}}{T},
\end{equation}
where $\gamma$ is a periodic orbit of $\phi_g^t$ (equivalently, a closed $g$-geodesic) and $l_g(\gamma)$ is its $g$-length. 
In particular, this formula says that in our setting, the total complexity of the flow's orbit structure can be read off of the complexity of the periodic orbit structure.

It is also important for the present paper that $h_{\rm top}(g)$ is a measure-theoretic entropy with respect to the \emph{measure of maximal entropy}. % (also known as the \emph{Bowen--Margulis measure}).
More precisely, 
for an invariant probability measure $\mu$, denote by $h_\mu(g)$ the measure-theoretic entropy with respect to $\mu$ of $\phi^1_g$. 
This quantity measures the exponential rate of growth (in $t$) of the number of all orbits of the geodesic flow distinguishable within time $t$ in a set of full $\mu$-measure. The two notions of entropy are related via the \emph{variational principle} \cite[Proposition 4.3.8]{FishHas}
 \begin{equation*}
\label{eq:BM}
\exists! \mu_{g}\in \mathcal P_{\phi^1_g}(S^gM),\quad h_{\mu_{g}}(\phi^1_g)=\sup_{\mu\in \mathcal P_{\phi_g^1}(S^gM)}h_{\mu}(\phi^1_g)=h_{\mathrm{top}}(\phi^1_g),
\end{equation*}
where $\mathcal P_{\phi^1_g}(S^gM)$ denotes the set of $\phi^1_g$-invariant probability measures on $S^gM$. %The 
%\emph{measure of maximal entropy} or 
The measure of maximal entropy of $(\phi^t_g)_{t\in \mathbb R}$ is the
unique invariant probability measure $\mu_{g}$ for which the supremum is attained.

The measure $\mu_g$ is a key object in this paper due to its presence in the formula for the first derivative of $h_{\rm top}(g)$ with respect to the metric. 
We first recall that by works of 
Katok, Knieper, Pollicott and Weiss, the map $\mathcal M^{\infty}_{<0}(M)\ni g\mapsto h_{\mathrm{top}}(g)$ is smooth \cite{KKW,KKPW,Pol}, where $\mathcal M^{\infty}_{<0}(M)$ is the space of $C^\infty$ Riemannian metrics of negative sectional curvature on $M$.
Moreover, using \eqref{eq:top}, Katok--Knieper--Weiss \cite{KKW} showed that if $(g_\la)_{\la\in (-\epsilon, \epsilon)}$ is a smooth path in $\mathcal M^{\infty}_{<0}(M)$, the first derivative is given by
\begin{equation}
\label{eq:KKW}\left.\frac{d}{d \la}\right|_{\la=0}h_{\mathrm{top}}(g_\la)=-\frac{h_{\mathrm{top}}(g_0)}2\int_{S^{g_0}M}[\partial_\la|_{\la=0} g_\la](v,v)d\mu_{g_0}(v).
\end{equation}
%Their proof uses Bowen's description of $\mu$ as the equidistribution of periodic orbits. 

\subsection{Symmetric tensors}
Recall that if $(g_{\lambda})_{\lambda \in (-\eps, \eps)}$ is a smooth path of metrics, then its derivative $\partial_{\lambda}|_{\lambda=0} g_{\lambda}$ is a symmetric two-tensor. As such, we now recall several standard notions about tensor fields; see, for instance, \cite[Chapter 14]{Lef} for further details.

Let $C^{\infty}(M;S^mT^*M)$ denote the smooth sections of the bundle of symmetric $m$-tensors on $M$. 
Note that the scalar product $g$ on $TM$ extends naturally to a scalar product on $C^{\infty}(M;S^mT^*M)$, which we will denote by $\langle \cdot, \cdot\rangle_{L^2(M;S^mT^*M)}$ (or $\langle \cdot,\cdot\rangle$ when there is no risk of confusion). It is defined as
\begin{equation}
\label{eq:PS}
    \forall S,T\in C^{\infty}(M;S^mT^*M), \ \langle S,T\rangle_{L^2(M;S^mT^*M)}:=\int_M \langle T_x,S_x\rangle d\mathrm{vol}_x(M), 
\end{equation}
where
$$\langle T_x,S_x\rangle=\sum_{1\leq i_1,\ldots, i_m\leq n}T(e_{i_1},\ldots,e_{i_m})S(e_{i_1},\ldots,e_{i_m}), $$
for any $g$-orthonormal basis $(e_i)_{i=1}^n$ of $T_xM$.
The \emph{trace} is given by
\begin{equation}
\label{eq:trace}
\tr_{g}: C^{\infty}(M;S^{m+2}T^*M)\to C^{\infty}(M;S^{m}T^*M),\quad S\mapsto \sum_{i=1}^nS_x(e_i,e_i,\cdot,\ldots,\cdot), 
\end{equation}
 where $(e_i)_{i=1}^n$ is a $g$-orthonormal basis of $T_xM$. 
 %add about tr(ST) = \langle T, S \rangle
 For $m=2$, the space of symmetric $2$-tensors decomposes as
\begin{equation}
\label{eq:decomp}C^{\infty}(M;S^2T^*M)=C^{\infty}(M;S^2_0T^*M)\oplus C^{\infty}(M)g, 
\end{equation}
where 
$
C^{\infty}(M;S^2_0T^*M):=\{S\in C^{\infty}(M;S^2T^*M)\mid \tr_{g}(S)=0\}
$
denotes the bundle of trace-free  symmetric $2$-tensors. 

Let $\nabla_g$ denote the Levi-Civita connection. We introduce the \emph{symmetrized covariant derivative}: 
$$D_{g}:=\mathrm{Sym}\circ \nabla_{g}: C^{\infty}(M;S^mT^*M)\to C^{\infty}(M;S^{m+1}T^*M).$$
The formal adjoint of $D_g$ is the \emph{divergence} operator $D_{g}^*$:
$$D_{g}^*=-\tr_g \circ \nabla_{g}: C^{\infty}(M;S^mT^*M)\to C^{\infty}(M;S^{m-1}T^*M). $$
The \emph{rough Laplacian} is given by
\begin{equation}
    \nabla^*_g\nabla_g: C^{\infty}(M;S^mT^*M)\to C^{\infty}(M;S^{m}T^*M).
\end{equation}
%{\color{blue}shoudl we mention it's the usual Laplacian on functions/conformal tensors? we use that a lot. we could also just note it commutes with the trace}
When there is no risk of confusion on the metric, we will suppress the $g$ subscripts.
\subsubsection{Spherical harmonics}\label{sec:spher-harm}
The relation between symmetric tensors and spherical harmonics is a powerful computational tool that we leverage in the proof of Proposition \ref{propHess1}.

The vertical Laplacian is defined, for any $f\in C^{\infty}(S^gM)$ by the formula 
$$\Delta_{\mathbb V}f(v)=\Delta_{S_x^gM}(f|_{S_x^gM}), \ \forall v\in T_xM, $$
where $\Delta_{S_x^gM}$ is the standard Laplacian on the fiber $S_x^gM\cong \mathbb S^{n-1}. $
We define the fiber bundle 
$$\Omega_m\longrightarrow M,\quad \Omega_m(x)=\mathrm{Ker}\big(\Delta_{S_x^gM}-m(n+m-2)\big), \ m\geq 0. $$
The sections of $\Omega_m$ are called \emph{spherical harmonics of degree $m$.}
Each smooth function $f\in C^{\infty}(S^gM)$ decomposes uniquely into a sum of spherical harmonics
$$
f=\sum_{m=0}^{+\infty}f_m, \quad  f_m\in \Omega_m \iff \Delta_{\mathbb V}f_m=m(n+m-2)f_m.
$$
 The relation with symmetric tensors is given by the pullback map: 
\begin{equation}
\label{eq:pim} \pi_m^*:C^{\infty}(M;S^mT^*M)\to C^{\infty}(SM), \quad \pi_m^*S(x,v)=S_x(\underbrace{v,\ldots,v}_{m \text{ times.}}).
\end{equation}
It defines an isomorphism from trace-free tensors of degree $m$ to spherical harmonics of degree $m$:
\begin{equation}
\label{eq:iso}
\pi_m^*: C^{\infty}(M,S_0^mT^*M)\overset{\sim}{\longrightarrow}  \Omega_m.
\end{equation}
We will need the following identity (see for instance \cite[Section 2.2]{GuKnLef}):
\begin{equation}\label{eq:pi-trace}
 \forall S\in C^{\infty}(M;S^2T^*M),\quad   \int_{SM} \pi_2^* S \, dm_g = \frac{1}{n {\rm Vol}_g(M)}\int_M {\rm tr}_g(S) \, d {\rm vol}_g,
\end{equation}
where $m_g$ is the \emph{Liouville measure} associated to $g$, i.e., the measure which is given locally by $d\mathrm{vol}_g \times d\mathrm{Leb}_{\mathbb S^{n-1}}$  and normalized so that we have a probability measure.

Any function $f\in C^{\infty}(M)$ belongs to $\Omega_0\subset C^{\infty}(S^gM)$. We will need the following relation between the $L^2(M, \mathrm{vol}_g)$-norm of $f$ and its $L^2(S^gM, m_g)$-norm,
\begin{equation}
    \label{eq:norms}
    \forall f,f'\in C^{\infty}(M),\quad \frac{1}{\mathrm{Vol}_g(M)}\langle f,f'\rangle_{L^2(M,\mathrm{vol}_g)}=\langle \pi_0^*f,\pi_0^*f'\rangle_{L^2(S^gM,m_g)}.
\end{equation}
\subsubsection{Action of the generator of the flow}
\label{sec:X}
We note the following relation between the symmetrized covariant derivative and the generator of the geodesic flow (see \cite[Lemma 14.1.9]{Lef}):
\begin{equation}
    \label{eq:XandD}
    X_g\pi_m^*=\pi_{m+1}^*D_g.
\end{equation}
The vector field $X_g$, seen as a derivation, is skew adjoint on $L^2(S^gM,dm_g)$, i.e., $X_g^*=-X_g.$  Moreover, for any $k\geq 0$:
$$X_g:C^{\infty}(M,\Omega_k)\to C^{\infty}(M,\Omega_{k-1})\oplus C^{\infty}(M,\Omega_{k+1}). $$
We will write $X_g=(X_g)_++(X_g)_-$, where $(X_g)_{\pm}:C^{\infty}(M,\Omega_k)\to C^{\infty}(M,\Omega_{k\pm1})$, see \cite[Lemma 14.2.1]{Lef}. We recall that $(X_g)_+^*=-(X_g)_-$.

\subsubsection{Decomposition of the space of symmetric tensors.}\label{decomp-slice}
 There exists a natural gauge given by the action of the group $\mathrm{Diff}_0(M)$ of smooth diffeomorphisms homotopic to the identity. We define
$$\mathcal O(g):=\{\phi^*g\mid \phi \in \mathrm{Diff}_0(M)\},\quad T_{g}\mathcal O(g)=\{\mathcal L_Vg\mid V\in C^{\infty}(M;TM)\}. $$
Note that any isometry invariant of $g$ is constant on orbits of the action, such as $h_{\rm top}(g)$ and the curvature tensors recalled in Section \ref{sec:curv} and appearing in the definition of the NRF.
Hence, we will  prove an injectivity result of the Hessian of the functional $\Phi$ defined in \eqref{eq:Phi} on a ``transverse slice" to $T_{g}\mathcal O(g)$. 
We note that
$$T_{g}\mathcal O(g)=\{D_{g}p\mid p\in C^{\infty}(M;T^*M)\}.$$

In particular, a natural transverse slice is provided by the kernel of the adjoint $D_{g}^*$ \cite[Theorem 14.1.11]{Lef}. Elements of $C^{\infty}(M;S^mT^*M)\cap \mathrm{Ker}(D_{g}^*)$ are called \emph{divergence-free} (or \emph{solenoidal}) tensors.
For any $S\in C^{\infty}(M;S^mT^*M)$, there exists a unique pair 
$$(p,h)\in C^{\infty}(M;S^{m-1}T^*M)\times \big(C^{\infty}(M;S^mT^*M)\cap \mathrm{Ker}(D_{g}^*)\big),\quad S=D_{g}p+h.$$

We recall the following lemma which allows one to ``project" a nearby metric $g'$ onto solenoidal tensors. It was obtained in this form in \cite[Lemma 2.4]{GuKnLef}, but the idea goes back to Ebin \cite{Eb}.
\begin{lemma}[Slice lemma]
\label{slice lemma}
Let $k\geq 2$,  and $\alpha\in (0,1)$. Then there exists a neighborhood $\mathcal U$ of $g$ in the $C^{k,\alpha}$-topology such that for any $g'\in \mathcal U$, there is a unique $\phi_{g '}\in \mathrm{Diff}_0(M)$ of regularity $C^{k+1,\alpha}$, close to identity, such that $\phi_{g'}^*g'\in \mathrm{Ker}(D_{g}^*)$ is divergence-free. Moreover, there exists $\epsilon>0$ and $C>0$ such that
$$\|g'-g\|_{C^{k,\alpha}}\leq \epsilon \ \Rightarrow \ \|\phi_{g'}^*g'-g\|_{C^{k,\alpha}}\leq C\|g'-g\|_{C^{k,\alpha}}.$$
\end{lemma}
\subsection{Curvature tensors}
\label{sec:curv}
We now recall the geometric notions appearing in the definition of the normalized Ricci flow.
 The \emph{Ricci tensor} is
\begin{equation*}
    \label{eq:Ric}
    \mathrm{Ric}_g\in C^{\infty}(M;S^2T^*M),\quad \mathrm{Ric}_g(v,w):=\mathrm{tr}(y\mapsto \mathbf{R}_g(v,y)w)=\sum_{i=1}^ng(\mathbf{R}_g(v,e_i)w,e_i),
\end{equation*}
for any $g$-orthonormal basis $(e_i)_{i=1}^n$, where $\mathbf{R}_g$ is the Riemannian curvature tensor of $g$. For $g_0$ hyperbolic and curvature $-1$, one has $\mathrm{Ric}_{g_0}=-(n-1)g_0.$ By a small abuse of notation, we will also denote by $\mathrm{Ric}_g$ the function $\pi_2^*(\mathrm{Ric}_g)\in C^{\infty}(S^gM)$.
The \emph{scalar curvature} is
\begin{equation*}
    \label{eq:Scalar}
    \mathrm{Scal}_g\in C^{\infty}(M),\quad \mathrm{Scal}_g(x):=\mathrm{tr}_g(\mathrm{Ric}_g)(x)=\sum_{i=1}^n\mathrm{Ric}_g(e_i,e_i),
\end{equation*}
for any $g$-orthonormal basis $(e_i)_{i=1}^n.$ For $g_0$ hyperbolic of curvature $-1$, one has $\mathrm{Scal}_{g_0}=-n(n-1).$ The \emph{total scalar curvature} is
\begin{equation}
    \label{eq:Scal}
    \mathcal S(g):=\int_M \mathrm{Scal}_{g}d\mathrm{vol}_g.
\end{equation}
%The \emph{normalized Ricci flow} (NRF) $(g_t)_{t\geq 0}$ starting from $g$ is defined by the differential equation
%\begin{equation}
%\label{eq:NRF-intro}
%al}{\partial{t}}g_t=-2\mathrm{Ric}_t +\frac{2}{n}\frac{\int_M\mathrm{Scal}_t d\mathrm{vol}_t }{\mathrm{Vol}_g(M)}g_t, \quad g_0=g,
%\end{equation} 
%where $\mathrm{Ric}_t$ is the Ricci tensor, $\mathrm{Scal}_t$ is the scalar curvature, $\mathrm{vol}_t$ is the volume form associated to $g_t$, and $\mathrm{Vol}_g(M)$ is the total volume of $M$ with respect to $g$. The average scalar curvature term on the right-hand side of \eqref{eq:NRF-intro} ensures the volume of $t \mapsto g_t$ remains constant and is the reason for the use of the word ``normalized".
Finally, we introduce the \emph{Lichnerowicz Laplacian}, which is given by
\begin{equation}\label{eq:lich}
(\Delta_L)_gS:=\nabla^*_g\nabla_gS+\mathrm{Ric}_g\circ S+S\circ \mathrm{Ric}_g-2R_g^\circ(S),
\end{equation}
where, for any $g$-orthonormal basis $(e_i)_{i=1}^n$, we have
\begin{equation*}
\label{eq:Ricci} R_g^{\circ}(S)(X,Y):=-\sum_{i=1}^n S(\mathbf{R}_g(e_i,X)Y,e_i),\quad S\circ \mathrm{Ric}_g(X,Y):=\sum_{i=1}^n S(\mathbf{R}_g(e_i,X)e_i, Y).
\end{equation*}
Note that when $g=g_0$ is a hyperbolic metric, one has
\begin{equation}
    \label{eq:LaplFla}
    \Delta_LS=\nabla^*\nabla S-2nS+2\tr(S)g_0,
\end{equation}
see \cite[Proof of Proposition 1.3.3]{Fla}.
The Lichnerowicz Laplacian appears naturally in our problem due to its relation to the variation of the Ricci curvature with respect to the metric (see, for instance, \cite{Bes}[Proposition 1.174 d]): 
\begin{equation}\label{eq:derRic}
      \partial_\la|_{\la=0} \mathrm{Ric}_\la=\frac 12(\Delta_L S)-D_{g_{0}}D_{g_{0}}^*S-\frac{1}{2}D_{g_0}d(\tr_{g_{0}}S).
\end{equation} 
%\subsection{The coercive estimate}
%\label{Xray}
%We recall a coercive estimate \cite{GL19} on the variance. 
%Let $f_1,f_2$ have zero mean, i.e., $\int_{SM}f_i(z)dm_g=0$ for $i=1,2$. The covariance of $f_1$ and $f_2$ with respect to the Liouville measure $m_g$ is
%\begin{equation}
%\label{eq:variance}
%\mathrm{Cov}_{g}(f_1,f_2):=\lim_{T\to +\infty}\frac 1 T\int_{SM}\left(\int_0^Tf_1(\phi_g^t(z))dt\right)\left(\int_0^Tf_2(\phi_g^t(z))dt\right)dm_g(z).
%\end{equation}
%When $f_1=f_2$, the covariance is called the variance and is written as $\mathrm{Cov}_{g}(f,f)=\mathrm{Var}_g(f).$
%In negative sectional curvature, we have the following coercive estimate due to Guillarmou and Lefeuvre, see for instance \cite[Lemma 2.1]{GuKnLef}. There is $C>0$, such that
%\begin{equation}
%\label{eq:coercive}
%\ \forall h\in H^{-1/2}(M;S^2T^*M)\text{ with } \int_{SM}\tr_g(h)dm_g=0,\quad \mathrm{Var}_g(\pi_2^*h) \geq C\|\Pi_{\mathrm{Ker}(D_{g}^*)}h\|_{H^{-1/2}}^2,
%\end{equation}
%where $\Pi_{\mathrm{Ker}(D_{g}^*)}$ is the orthogonal projection on solenoidal (or divergence-free) tensors.

\section{Definition and first properties of the functional $\Phi$}
\label{sec:1stproperties}
In this section, we consider a smooth negatively curved metric $g$, together with the normalized Ricci flow $(g_t)_{t\geq 0}$ starting from $g$ (see \eqref{eq:NRF-intro}). We further suppose that $g$ is $C^2$-close to a hyperbolic metric $g_{\rm hyp}$. By work of Ye \cite{Ye}, this implies $(g_t)_{t\geq 0}$ is well-defined and negatively curved for all $t\geq 0$, and converges to a hyperbolic metric as $t\to+\infty$.
\subsection{Smoothness of the derivative}
From the results recalled in Section \ref{section: topological entropy}, we know that $t\mapsto h_{\rm top}(g_t)$ is smooth. From \eqref{eq:KKW} and \eqref{eq:NRF-intro}, we get 
\begin{align}\label{derivative of top}
\frac{d}{dt}h_{\mathrm{top}}(g_t) &= -\frac{h_{\mathrm{top}}(g_t)}{2}\int\limits_{S^{g_t}M}\pi_2^*\partial_t g_t\,d\mu_{g_t}\nonumber\\
&= -\frac{h_{\mathrm{top}}(g_t)}{2}\left(\int_{S^{g_t}M} -2{\rm {Ric}}_{g_t} (v) \, d \mu_{g_t} + \frac{2}{n {\rm {Vol}}_{g_t}(M)} \mathcal{S}(g_t)\right).
\end{align}
We define a functional $\Phi$ by
\begin{equation}
\label{eq:Phi}
\begin{split}
\Phi(g) &:=  \int_{S^{g}M} -{\rm {Ric}}_g (v) \, d \mu_g + \frac{1}{n {\rm {Vol}}_g(M)} \mathcal{S}(g).
\end{split}
\end{equation}
Note that if $t \mapsto \Phi(g_t)$ is positive, then $t \mapsto h_{\rm top}(g_t)$ is (strictly) decreasing. In order to show that $\Phi$ is positive in a neighborhood of $g_{\rm hyp}$, we will show in Section \ref{sub:crit} that $g_{\rm hyp}$ is a critical point of $\Phi$ and then compute the Hessian of $\Phi$ at $g_{\rm hyp}$ in Section \ref{sec:Hessian}. Before taking derivatives of $\Phi$, we justify that $\Phi$ is smooth in $g$ (Corollary \ref{cor:Phi-smooth}).

To carry out these analyses, we identify the different unit tangent bundles $S^gM$ with $S^{g_{\rm hyp}} M$ by rescaling each fiber:
\begin{equation}
    \label{eq:Phi_g}
    \Psi_g:SM:=S^{g_{\rm hyp}}M\to S^gM,\quad (x,v)\mapsto \big(x,\tfrac{v}{\|v\|_g} \big).
\end{equation}
Using $\Psi_g$, we pull back the geodesic vector field and the measure of maximal entropy of $g$ to $SM$:
\begin{equation}
    \label{eq:pb}
    \tilde{X}_g:=\Psi_g^*X_g\in C^{\infty}(SM;T(SM)), \quad \tilde \mu_g:= \Psi_g^*\mu_{g_{\rm}}.
    \end{equation}
    For any metric $g$, we have that $\tilde \mu_g$ is a $\tilde X_g$-invariant probability measure on $SM.$ In this notation, using that $\mathrm{Ric}_g$ is a $2$-tensor, we obtain
\begin{equation}
\label{eq:Phi2}
    \Phi(g)=\int_{SM} \frac{1}{\Vert v \Vert_{g}^2}(- {\rm Ric}_g (v)) \, d \tilde \mu_g + \frac{1}{n {\rm Vol}_g(M)} \mathcal{S}(g).
\end{equation}
We also record for future use the following lemma which computes the variation of the pullback 
$\tilde m_g:=\Psi_g^*m_g$ of the Liouville measure $m_g$ of $g$ to $SM$.
\begin{lemma}\label{lem:liou} Let $(g_\la)_{\la\in (-\epsilon,\epsilon)}$ be a perturbation of a metric. Then the derivative of the $($normalized$)$ pulled back Liouville measures is
     $$\partial_{\lambda}|_{\la=0} d \tilde m_{\lambda} = \big( {\rm tr} (\partial_\la|_{\la=0} g_{\lambda})-\frac n 2\pi_2^*(\partial_\la|_{\la=0} g_{\lambda})-\tfrac 1{2\mathrm{Vol}_{0}(M)}\int_{M}\tr_{0}(\partial_\la|_{\la=0} g_{\lambda})d\mathrm{vol}_{0}\big) \, d  m_0,$$ where the subscript $\la$ denotes the objects corresponding to $g_\la.$
\end{lemma}

\begin{proof}
Let $\Omega_{\la}$ be the Sasaki volume form of $g_{\la}$ on $S^{g_\la} M$. Recall that the Liouville measure $m_{\la}$ is the (probability) measure induced by the volume form $\frac{1}{{\rm Vol}(g_{\la})}\Omega_{\la}$.
Differentiating the pullback $\frac{1}{{\rm Vol}(g_{\la})} \Psi_{\la}^* \Omega_{\la}$ at $\la = 0$ and using the expression for $\Psi_{g_\la}^* \Omega_{\la}$ obtained at the end of the proof of \cite[Lemma 6.4]{CL25} yields the desired result.
%The variation of pullback Liouville measure was computed in \cite[Lemma 6.4]{CL25}.
%More precisely, we obtain the formula stated in the lemma by differentiating the expression of the pullback Liouville measure in local coordinates obtained in the end of the proof of \cite[Lemma 6.4]{CL25}
\end{proof}    

To check the smoothness of $\Phi$, we use equation \eqref{eq:Phi2}. 
The following lemma shows the smoothness of the (pulled back)
%or pullback *of the* BM measure, if you prefer that wording
Bowen-Margulis measure $\tilde \mu_ g$ with respect to the metric~$g$. Throughout the rest of the paper, we will denote by $\mathcal M_{<0}^k(M)$ the set of $C^k$ Riemannian metrics of negative curvature.
\begin{lemma}[cf. Theorem C c) in \cite{Con}]\label{lemm:BMsmooth} Let $k\geq 0$ and let $f\in C^{\infty}(SM)$ be fixed. Then
\begin{equation}
\label{eq:musmooth}
  \mathcal M_{<0}^{k+2}(M)\ni g   \mapsto \int_{SM} f(x)d\tilde \mu_g (x),
\end{equation}
is a $C^k$ map.
\end{lemma}

Using Lemma \ref{lemm:BMsmooth} and \eqref{eq:Phi2}, we deduce the following
\begin{corollary}\label{cor:Phi-smooth}
$\mathcal M^{\infty}_{<0}(M)\ni g\mapsto  \Phi(g)$ is smooth.
\end{corollary}

\begin{remark} 
We will use a more refined version of Lemma \ref{lemm:BMsmooth} (i.e., \cite[Theorem C c)]{Con}) to determine the error term in our second-order Taylor expansion of $\Phi$ in Lemma 5.1; see \eqref{eq:contreras-map}. 
We note, however, that the result of Corollary \ref{cor:Phi-smooth} can be deduced in several other ways. Firstly, it was previously known by \cite[Corollary 1]{KKPW} that the map \eqref{eq:musmooth} is $C^{k-1}$ when $k\geq 1.$
Also, one can show that \eqref{eq:musmooth} is actually a $C^{k+1}$ map using the theory of \emph{Pollicott--Ruelle} resonances. 
%We omitted the detailed proof as it gives less explicit bounds on the regularity needed (in the variation of metrics) to bound the derivatives of $\Phi$ in Lemma~\ref{lemm:remainder} in comparison to \cite[Theorem C c)]{Con}. 
Indeed, it is shown in \cite[Theorem 1]{Hum1} that $h_{\mathrm{top}}(g)$ is a simple\footnote{Here, we mean that it has geometric multiplicity equal to $1$ and no Jordan block.} eigenvalue (called a Pollicott--Ruelle resonance) for the action of $\tilde X_g$ on $(n-1)$-forms. Moreover, the measure $\tilde \mu_g$ can be obtained as the distributional product of the corresponding eigenvector (called the \emph{resonant state}) and the dual eigenvector (the \emph{coresonant state}); see \cite[Theorem 2]{Hum1}. Since the eigenvalue is simple, $\tilde \mu_g$ is equal to the trace of the spectral projector at $h_{\mathrm{top}}(g)$. The spectral projector varies smoothly in $g$ (with no loss of derivative), and we deduce the desired smoothness of $g\mapsto \tilde \mu_g.$ (We omit the detailed version of the above argument as it gives less explicit bounds on the regularity needed (in the variation of metrics) to bound the derivatives of $\Phi$ in Lemma~\ref{lemm:remainder} in comparison to \cite[Theorem C c)]{Con}.)
\end{remark}

\subsection{Hyperbolic metrics are critical points of $\Phi$}\label{sub:crit}
Since hyperbolic metrics are fixed points of the normalized Ricci flow, the functional $\Phi$ vanishes at any hyperbolic metric. Moreover, we show that $\Phi$ has a critical point at any hyperbolic metric:

\begin{lemma}
\label{lemm:crit}
    Let $(M^n,g_{\rm hyp})$ be a closed hyperbolic manifold and let $(g_{\la})_{\la \in (-\eps, \eps)}$ be a smooth perturbation of $g_0: = g_{\rm hyp}$. Then $\partial_{\la}|_{\la = 0} \Phi(g_{\la}) = 0$. 
\end{lemma}
\begin{proof} Write $S:=\partial_\la|_{\la=0}g_\la.$ Differentiating \eqref{eq:Phi2} gives
\begin{align*}
\partial_\la|_{\la=0}\Phi(g_\la)
&=(n-1)\int_{SM}\partial_\la|_{\la=0}\left(\frac{1}{\Vert v \Vert_{\la}^2}\right)  \, d m_0
-\int_{SM}\partial_\la|_{\la=0}{\rm Ric}_\la \,dm_0 \\
&+(n-1)\int_{SM}\partial_\la|_{\la=0}\,d\tilde \mu_\la -\frac{1}{n}\int_M n(n-1) \partial_{\la}|_{\la = 0} \left( \frac{d {\rm vol}_{\la}}{{\rm Vol}_{\la}(M)} \right)\\
&+  \frac{1}{n\mathrm{Vol}(M)}\mathcal \int_M \partial_\la|_{\la=0} {\rm Scal}_{\la} d {\rm vol}_{0}\\
\\&=-(n-1)\int_{SM}\pi_2^*S(v)\,dm_0(v) +  \frac{1}{n {\rm Vol}(M)} \int_{M}  \partial_{\la}|_{\la = 0} {\rm tr}_{\la} (-(n-1) g_0)\,d{\rm vol}_0.
\end{align*}
where in the for last equality we used \eqref{eq:pi-trace}, the fact that $\tilde \mu_\la$ and $\frac{1}{{\rm Vol}_{\la}(M)}  {\rm vol}_{\la}$ are probability measures (on $SM$ and $M$, respectively), that ${\rm Scal}_g = {\rm tr}_g {\rm Ric}_g$ and ${\rm Ric}_{g_0}=-(n-1)g_0.$

To finish the proof, we recall that for any $2$-tensor $T$ and any $x\in M$ we have (see, for example, \cite[Proposition 2.3.6]{Topping})
$$\partial_\lambda|_{\lambda=0}\tr_\lambda(T_x)=-\langle S,T\rangle_x,$$
where the right-hand side is defined in \eqref{eq:PS}.
Using this followed by \eqref{eq:pi-trace} on the second term above shows that $\partial_{\la}|_{\la = 0} \Phi(g_{\la}) = 0$ for any perturbation of $g_0 := g_{\rm hyp}$.
\end{proof}

\section{Computing the Hessian of $\Phi$ at a hyperbolic metric}
\label{sec:Hessian}
In this section, we compute the second variation of $\Phi$ at $g_{\rm hyp}$.
%{\color{blue} added the following:}
By Lemma \ref{lemm:crit}, $g_{\rm hyp}$ is a critical point of $\Phi$, which means that for any perturbation $(g_{\la})_{\la \in (-\eps, \eps)}$ of $g_0 := g_{\rm hyp}$, 
we have that $\partial^2_{\la}|_{\la = 0} \Phi(g_{\la})$ depends only on 
$S: = \partial_{\la}|_{\la = 0} \, g_{\la}$ (and not on $\partial^2_{\la}|_{\la = 0}g_{\la}$). More precisely, we have $\partial^2_{\la}|_{\la = 0} \Phi(g_{\la}) = d^2_{g_{\rm hyp}} \Phi(S, S)$.

In order to compute $\partial^2_{\la}|_{\la = 0} \Phi(g_{\la})$, we may (and will) always assume that $g_{\la}$ is a \emph{first-order perturbation}, i.e., it is of the form $g_{\la} = g_0 + \lambda S$ for some symmetric 2-tensor $S$. In this section, we compute $\partial^2_{\la}|_{\la = 0} \Phi(g_{\la})$ and show that it is positive for $S$ divergence-free (defined in Section \ref{decomp-slice}) and of \emph{zero mean}, 
that is, $\int_{SM} \pi_2^*S \, dm_0 =0$. 
(Note that by \eqref{eq:pi-trace} and \cite[Proposition 1.186 b)]{Bes}, the zero-mean condition is equivalent to $\partial_{\la}|_{\la = 0} {\rm Vol}_{\la}(M) = 0$.)

\subsection{Computing the Hessian}
First, we introduce the following notation. Let $S=\partial_\la|_{\la=0}g_\la$ where $(g_\la)_{\la\in (-\epsilon,\epsilon)}$ is a deformation of $g_0$. 
We let 
\begin{equation}
    \label{eq:dmu}
    \mathcal D_{g_0}\tilde \mu(S):=\partial_\la|_{\la=0}d\tilde \mu_\la,
\end{equation}
which is well defined by Lemma \ref{lemm:BMsmooth}. In the following, we will denote by $(f,\mathcal D_{g_0}\tilde \mu(S))$ the distributional pairing with any smooth function $f\in C^\infty(SM).$

In this section, we prove the following proposition.
\begin{proposition}
\label{propHess1}
    Let $(M^n,g_{\rm hyp}) $ be a closed hyperbolic manifold, let $S$ be a divergence-free tensor of zero mean, 
    and let $(g_\la)_{\la \in (-\eps,\eps)}$ be a variation such that $S=\partial_\la|_{\la=0}g_\la$. Then
   \begin{align*}
\partial^2_\la|_{\la=0}&\Phi(g_\la) = 2(n-1)\left(\|\pi_2^*S\|^2-\frac{1}{n\mathrm{Vol}(M)}\big(\| S\|^2 -\frac{1}{2}\|\tr_{g_0}(S)\|^2\big)-(\pi_2^*S(v), \mathcal D_{g_0}\tilde \mu(S))\right)
\\+&\left(\langle \pi_2^*\Delta_LS,\pi_2^*S\rangle -\frac{1}{n\mathrm{Vol}(M)}\big(\langle \Delta_L S, S\rangle -\frac{1}{2}\langle \tr_{g_0}(\Delta_L S),\tr_{g_0}(S)\rangle\big)-(\pi_2^*(\Delta_L S),\mathcal D_{g_0}\tilde \mu(S))\right).\end{align*}
%where the subscript $0$ denotes the objects corresponding to $g_0=g_{\rm hyp}.$
\end{proposition}

\begin{proof}
As explained at the beginning of this section, we may assume that $g_{\la} = g_{\rm hyp}+ \lambda S$ for $S$ as in the statement of the proposition.  
We compute the second derivative of \eqref{eq:Phi2}, using that $\mathrm{Ric}_0=-(n-1)g_0$ and $m_0=\mu_0$. Everywhere below the subscript $\lambda$ denotes the objects corresponding to $g_\lambda$. We have
\begin{align*}
\partial^2_\la|_{\la=0}&\Phi(g_\lambda) = \partial^2_\lambda|_{\la=0}\left(\int_{SM} \frac{1}{\Vert v \Vert_{\la}^2}(- {\rm Ric}_\la (v)) \, d \tilde \mu_\la + \frac{1}{n {\rm Vol}_\la(M)} \mathcal{S}(\la)\right)
\\&= (n-1)\int_{SM} \partial^2_\lambda|_{\la=0}\left(\frac{1}{\Vert v \Vert_{\la}^2}\right) \, d m_0 
+ (n-1)\int_{SM}  \,\partial^2_\la|_{\la=0} d \tilde \mu_\la 
\\&+\int_{SM} \partial_{\la}^2|_{\la = 0} (- {\rm Ric}_\la (v)) \, d m_0 + \partial_{\lambda}^2 |_{\lambda =0} \left(\frac{1}{n {\rm Vol}_{\la}(M)} \mathcal{S}(\lambda)\right)
\\& +2(n-1)\left( \partial_{\la} |_{\la = 0}\left( \frac{1}{\Vert v \Vert_{\la}^2} \right) , \mathcal D_{g_0}\tilde \mu(S) \right)
+ 2 \Big( \partial_{\la}|_{\la = 0} (- {\rm Ric}_\la (v)) , \mathcal D_{g_0}\tilde \mu(S) \Big) \\
&+2\int_{SM} \partial_{\la}|_{\la = 0}\left( \frac{1}{\Vert v \Vert_{\la}^2} \right) \partial_{\la}|_{\la = 0}(- {\rm Ric}_\la (v)) \, d m_0=:A+B+C+D.
\end{align*}
We compute each of $A,B,C,D$ separately. First, in order to compute $A$, note that we have
\begin{equation}\label{eq: derivatives of the inverse norm of v}
 \partial_{\la}|_{\la=0}\left(\frac{1}{\|v\|_\la^2}\right)=- \pi_2^*S,\quad \partial_{\la}^2|_{\la=0}\left(\frac{1}{\|v\|_\la^2}\right)=-\pi_2^*({\partial_{\la}^2|_{\la=0}g_\la})
 +2(\pi_2^*S)^2 = 2(\pi_2^*S)^2,
 \end{equation}
 where we used that $\partial_{\la}|_{\la = 0}^2g_\la=0$ since the perturbation is linear. 
 Everywhere below $\tr$ is $\tr_0$ if applied to a symmetric $2$-tensor. Using \eqref{eq: derivatives of the inverse norm of v}, \eqref{eq:pi-trace} and the fact that $\tilde \mu_\la$ is a probability measure for any $\la$, we have
\begin{equation}
\label{eq:A}
\begin{split}
    A&=(n-1)\int_{SM} \partial^2_\lambda|_{\la=0}\left(\frac{1}{\Vert v \Vert_{\la}^2}\right) \, d m_0 
+ (n-1)\int_{SM}  \,\partial^2_\la|_{\la=0} d \tilde \mu_\la
\\&=2(n-1)\|\pi_2^*S\|^2+(n-1)\partial^2_\la|_{\la=0}\tilde \mu_\la(SM)
\\&=2(n-1)\|\pi_2^*S\|^2.
\end{split}
\end{equation}
Next, we compute $B$. Using \eqref{eq:pi-trace} and \eqref{eq:Scalar} and recalling that $\mathrm{Ric}_0=-(n-1)g_0$ and $\mathrm{Scal}_0=-n(n-1)$, we obtain
\begin{align*}
    B&=\int_{SM} \partial_{\la}^2|_{\la = 0} (- {\rm Ric}_\la (v)) \, d m_0 + \partial_{\lambda}^2 |_{\lambda =0} \left(\frac{1}{n {\rm Vol}_{\la}(M)} \mathcal{S}(\lambda)\right)
    \\&=\frac{-(n-1)}{n {\rm Vol}(M)}\int_M \partial_\la^2|_{\la=0}\tr_\la (g_0)\,d\mathrm{vol}_0 -(n-1)\int_{M}\partial^2_\la|_{\la=0}\left(\frac{d\mathrm{vol}_\la}{{\rm Vol_{\la}(M)}}\right)
    \\&-\frac{2(n-1)}{n {\rm Vol}(M)}\int_M \partial_\la|_{\la=0} \tr_\la (g_0)\,\partial_\la|_{\la=0}d\mathrm{vol}_\la+\frac{2}{n {\rm Vol}(M)}\int_M  \tr(\partial_\la|_{\la=0}\mathrm{Ric}_\la)\,\partial_\la|_{\la=0}d\mathrm{vol}_\la
    \\&+\frac{2}{n {\rm Vol}(M)}\int_M \partial_\la|_{\la=0} \tr_\la ( \partial_\la|_{\la=0}\mathrm{Ric}_\la) \,d\mathrm{vol}_0,
\end{align*}
where we used the mean zero condition, i.e., that $\partial_{\la}|_{\la = 0}{\rm Vol}_{\la}(M) = 0$, to deduce that 
\begin{align*}
\frac{2}{n}\int_M\partial_\lambda|_{\lambda=0}\big(\tr_\lambda({\rm Ric}_\lambda(M))\big)&\partial_\lambda|_{\lambda=0}\left(\frac{d{\rm vol}_\lambda}{{\rm Vol}_\lambda(M)}\right)\\&=\frac{2}{n{\rm Vol}(M)}\int_M\partial_\lambda|_{\lambda=0}\big(\tr_\lambda({\rm Ric}_\lambda(M))\big)\partial_\lambda|_{\lambda=0}d{\rm vol}_\lambda.
\end{align*}
Since $\frac{\mathrm{vol}_{\la}}{\mathrm{Vol}_\la(M)}$ is a probability measure, one has $\int_{M}\partial^2_\la|_{\la=0}\left(\frac{d\mathrm{vol}_\la}{{\rm Vol}_{\la}(M)}\right)=0.$ %{\color{teal} changed the previous two lines to make the computation work without assuming constant volume.}

Next, we recall the following equalities for any $2$-tensor $T$ and any $x\in M$ (see, for example, \cite[Proposition 2.3.6]{Topping}):
$$\partial_\lambda|_{\lambda=0}\tr_\lambda(T_x)=-\langle S,T\rangle_x\qquad\text{and}\qquad \partial^2_\lambda|_{\lambda=0}\tr_\lambda(T_x)=2\langle S,S\circ T\rangle_x-\langle\ddot{g}_0,T\rangle_x=2\langle S,S\circ T\rangle_x,$$
where $S\circ T$ denotes the composition of $(2,0)$-tensors seen as $(1,1)$-tensors and where we used that $\ddot{g}_0=0$ for the linear perturbation $(g_\la)_{\la \in (-\eps,\eps)}$. Using this, together with \eqref{eq:PS}, the fact that $\partial_\la|_{\la=0}d\mathrm{vol}_\la=\tfrac 12\tr(S)d\mathrm{vol}_0$ (see, for instance, \cite[Proposition 1.186 a)]{Bes}) and \eqref{eq:derRic}, we obtain 
\begin{align*}
    B&=\frac{-2(n-1)}{n {\rm Vol}(M)}\|S\|^2+\frac{n-1}{n {\rm Vol}(M)}\|\tr(S)\|^2
    \\&+\frac{1}{2n {\rm Vol}(M)}\int_M  \tr(\Delta_L S-D_{g_0}d\tr(S))\tr(S)\,d\mathrm{vol}_0-\frac{1}{n {\rm Vol}(M)}\langle S,\Delta_L S-D_{g_0}d\tr(S)\rangle.
\end{align*}
Using \eqref{eq:norms}, we further simplify the above expression using the decomposition of symmetric tensors into spherical harmonics; see Sections \ref{sec:spher-harm} and \ref{sec:X}. More precisely, denoting by $\Pi_{\Omega_0}$ the orthogonal projection onto $\Omega_0$, and using \eqref{eq:norms}, together with \eqref{eq:XandD}, $X=X_++X_-$ and $(X_+)^*=-X_-$, we obtain
\begin{align*}
    &\frac{-1}{2n {\rm Vol}(M)}\int_M \tr(D_{g_0}d\tr(S))\tr(S)d\mathrm{vol}_0
=\frac{-1}{2n}\langle  \pi_0^* \tr(D_{g_0}d\tr(S)),\pi_0^*\tr(S)\rangle_{L^2(SM)}
\\&=-\frac 12\langle  \Pi_{\Omega_0} \pi_2^*D_{g_0}d\tr(S),\pi_0^*\tr(S)\rangle=-\frac 12\langle X_-X_+ \pi_0^*\tr(S),\pi_0^*\tr(S)\rangle_{L^2(SM)}=\frac 12\|X(\pi_0^*\tr(S))\|^2. \end{align*}
Next, using the fact that $S$ is divergence-free, we have
$$\frac{1}{n {\rm Vol}(M)}\langle S,D_{g_0}d\tr(S)\rangle=\frac{1}{n {\rm Vol}(M)}\langle D_{g_0}^*S,d\tr(S)\rangle=0. $$
In total, we have shown
\begin{equation}
    \label{eq:B}
    \begin{split}
    B&=\frac{-2(n-1)}{n {\rm Vol}(M)}\|S\|^2+\frac{n-1}{n {\rm Vol}(M)}\|\tr(S)\|^2
    \\&+\frac{1}{2n {\rm Vol}(M)}\langle \tr(\Delta_L S),\tr(S)\rangle -\frac{1}{n\mathrm{Vol}(M)}\langle S,\Delta_L S\rangle +\frac 1 2 \|X(\pi_0^*\tr(S))\|^2.
    \end{split}
\end{equation}
We now simplify $C$. Using \eqref{eq: derivatives of the inverse norm of v} and \eqref{eq:derRic}, 
\begin{align*}
    C&=2(n-1)\left( \partial_{\la} |_{\la = 0}\left( \frac{1}{\Vert v \Vert_{\la}^2} \right), \mathcal D_{g_0}\tilde \mu(S)\right)
+ 2 \big( \partial_{\la}|_{\la = 0} (- {\rm Ric}_\la (v)) , \mathcal D_{g_0}\tilde \mu(S)\big)
\\&=-2(n-1)(\pi_2^*S, \mathcal D_{g_0}\tilde \mu(S))-(\pi_2^*(\Delta_L S), \mathcal D_{g_0}\tilde \mu(S))+(\pi_2^*(D_{g_0}d\tr(S)), \mathcal D_{g_0}\tilde \mu(S)).
\end{align*}
To simplify the third term above, we note that for any function $f\in C^{\infty}(SM)$, the $\tilde X_\la$-invariance of $\tilde \mu_\la$ gives
$$0=\partial_{\la}|_{\la=0}\left(\int_{SM}\tilde X_\la(f)d\tilde \mu_\la\right)=( X(f),\mathcal D_{g_0}\tilde \mu(S))+\int_{SM} \partial_{\la}|_{\la=0}\tilde X_\la(f)d \mu_0. $$
 Using the $\tilde X_\la$-invariance of $\tilde m_\la$, the pullback of the Liouville measure $m_\la$ of $g_\la$, a similar computation gives
$$0=\partial_{\la}|_{\la=0}\left(\int_{SM}\tilde X_\la(f)d\tilde m_\la\right)=\int_{SM} X(f)\partial_\la|_{\la=0}d\tilde m_\la+\int_{SM} \partial_{\la}|_{\la=0}\tilde X_\la(f)d m_0.  $$ Since $dm_0=d\mu_0$, we obtain
\begin{equation}
    \label{eq:smart}
    (X(f),\mathcal D_{g_0}\tilde \mu(S))=-\int_{SM}\partial_\la|_{\la=0}\tilde X_\la(f)dm_0=\int_{SM} X(f)\partial_\la|_{\la=0}d\tilde m_\la.
\end{equation}
Using Lemma \ref{lem:liou} and that $S$ has mean zero, we have $\partial_\la|_{\la=0}d\tilde m_\la=( \tr(S)-\tfrac n2 \pi_2^*S)dm_0$.
Since $\pi_2^*(D_{g_0}d\tr(S))=X^2(\pi_0^*\tr(S))$, we can apply \eqref{eq:smart} to $f=X(\pi_0^*\tr(S))$ to obtain, using $X^*=-X,$
\begin{align*}
(\pi_2^*(D_{g_0}d\tr(S)), \mathcal D_{g_0}\tilde \mu(S))&=\int_{SM}X^2(\pi_0^*\tr(S))\pi_0^*\tr(S)dm_0-\frac n2\int_{SM} X^2(\pi_0^*\tr(S))\pi_2^*Sdm_0
\\&=-\|X(\pi_0^*\tr(S))\|^2+\frac{n}{2}\langle X(\pi_2^*S),X(\pi_0^*\tr(S))\rangle.
\end{align*}
We now simplify the last term. Using the fact that the spaces $\Omega_m$ are orthogonal for different values of $m$, we obtain
\begin{align*} 
\langle X(\pi_2^*S),X(\pi_0^*\tr(S))\rangle&=\langle X_-(\pi_2^*S),X_+(\pi_0^*\tr(S))\rangle+\langle X_+(\pi_2^*S),X_+(\pi_0^*\tr(S))\rangle
\\&=\langle X_-(\pi_2^*S_0),X_+(\pi_0^*\tr(S))\rangle+\langle X_+(\tfrac 1n \pi_0^*\tr(S)),X_+(\pi_0^*\tr(S))\rangle
\\&=\langle X_-(\pi_2^*S_0), X_+(\pi_0^*\tr(S))\rangle+\tfrac{1}{n}\|X_+(\pi_0^*\tr(S))\|^2,
\end{align*}
where $S_0\in C^{\infty}(M;S^2_0T^*M)$ is the trace-free part of $S.$
Since $S$ is divergence-free, we have (see, for instance, \cite[$(2.14)$]{Hum})
$$X_-(\pi_2^*S_0)=-\frac{2}{n(n+2)}X_+(\pi_0^*\tr(S)).$$
 This implies
\begin{align}
\label{eq:divfreeused}\langle X(\pi_2^*S),X(\pi_0^*\tr(S))\rangle&=\frac{1}{n}\left(1-\frac{2}{n+2}\right)\|X_+(\pi_0^*\tr(S))\|^2
\\&=\frac{1}{n+2}\|X_+(\pi_0^*\tr(S))\|^2 = \frac{1}{n+2}\|X(\pi_0^*\tr(S))\|^2.\nonumber
\end{align}
In total, we have shown that 
\begin{equation}
    \label{eq:C}
    \begin{split}
        C=-2(n-1)(\pi_2^*S, \mathcal D_{g_0}\tilde \mu(S))-(\pi_2^*(\Delta_L S), \mathcal D_{g_0}\tilde \mu(S))-\frac{n+4}{2(n+2)}\|X(\pi_0^*\tr(S))\|^2.
    \end{split}
\end{equation}
Finally, we simplify $D.$ We have
\begin{align*}
    D&=2\int_{SM} \partial_{\la}|_{\la = 0}\left( \frac{1}{\Vert v \Vert_{\la}^2} \right) \partial_{\la}|_{\la = 0}(- {\rm Ric}_\la (v))dm_0=\langle \pi_2^*S,\pi_2^*(\Delta_L S-D_{g_0}d\tr(S))\rangle
    \\&=\langle \pi_2^*S,\pi_2^*\Delta_L S\rangle -\langle \pi_2^*S,\pi_2^*(D_{g_0}d \tr(S))\rangle.
\end{align*}
Using $\pi_2^*(D_{g_0}d \tr(S))=X^2(\pi_0^*\tr(S))$, the fact that $X^*=-X$, the decomposition $X=X_++X_-$, we see that 
$$ -\langle \pi_2^*S,\pi_2^*(D_{g_0}d \tr(S))\rangle=-\langle \pi_2^*S, X^2(\pi_0^*\tr(S))\rangle=\langle X(\pi_2^*S),X(\pi_0^*\tr(S))\rangle.$$
Using, \eqref{eq:divfreeused}, we have shown that 
\begin{equation}
    \label{eq:D}
    D=\langle \pi_2^*S,\pi_2^*\Delta_L S\rangle+\frac{1}{n+2}\|X(\pi_0^*\tr(S))\|^2.
\end{equation}
Adding \eqref{eq:A}, \eqref{eq:B}, \eqref{eq:C}, \eqref{eq:D} completes the proof of the proposition. 
\end{proof}
\subsection{Computing the derivative of the Bowen--Margulis measure}
We now analyze the derivative of the Bowen--Margulis measure. The variation of the Bowen--Margulis measure with respect to the metric is a subtle question. Indeed, for a metric close (but not isometric) to $g_{\rm hyp}$, the Liouville and Bowen--Margulis measures are not equal and, therefore, mutually singular \cite{Fla,Hum}.

First, we recall the definition of the \emph{covariance}. Let $f_1,f_2$ have zero mean, meaning, $\int_{SM}f_i(z)dm_g=0$ for $i=1,2$. The covariance of $f_1$ and $f_2$ with respect to the Liouville measure $m_g$ is
\begin{equation}
\label{eq:variance}
\mathrm{Cov}_{g}(f_1,f_2):=\lim_{T\to +\infty}\frac 1 T\int_{SM}\left(\int_0^Tf_1(\phi_g^t(z))\,dt\right)\left(\int_0^Tf_2(\phi_g^t(z))\,dt\right)\,dm_g(z).
\end{equation}
When $f_1=f_2$, the covariance is called the \emph{variance} and is written as $\mathrm{Cov}_{g}(f,f)=:\mathrm{Var}_g(f).$
We say a symmetric $2$-tensor $S$ has \emph{zero mean} if $\int_{SM} \pi_2^*S\, dm_0 = 0$.
Using computations due to Flaminio \cite{Fla} (see also Pollicott \cite{Pol}), we show the following result. 
\begin{proposition}\label{lemma:Flaminio}
Let $S_1,S_2$ be two symmetric $2$-tensors with zero mean. 

Then we have
\begin{align}\label{eq:polarization}
    \frac 12\big(&{\rm Cov}(\pi_2^* T(S_1), \pi_2^* S_2)+{\rm Cov}(\pi_2^* T(S_2), \pi_2^* S_1)\big)
    \\&= \langle \pi_2^* S_1,  \pi_2^* S_2 \rangle 
    - \frac{1}{n {\rm Vol}(M)} (\langle S_1, S_2 \rangle - \frac{1}{2} \langle {\rm tr}(S_1), {\rm tr}(S_2) \rangle )- (\pi_2^* S_1, \mathcal{D}_{g_0} \tilde \mu(S_2)),\nonumber
\end{align}    
    where $T(S) = - \frac{1}{2} S + \frac{1}{4} \nabla^* \nabla S + \frac{1}{2} {\tr}_{g_0}(S)g_0-\frac{1}{2}D_{g_0}D_{g_0}^*S$.
\end{proposition}

\begin{proof} 
The first step of the proof is to show that all terms in \eqref{eq:polarization} are bilinear and symmetric.
Then, we will show that \eqref{eq:polarization} holds for $S_1= S_2$. Finally, we will deduce the general form of \eqref{eq:polarization} using the polarization identity.

For the first step, we need to show that $(S_1, S_2) \mapsto (\pi_2^* S_1, \mathcal{D}_{g_0} \tilde \mu(S_2))$ is bilinear and symmetric. To see this,
    let $g_{\la_1, \la_2}$ be any two-parameter family satisfying
\begin{equation}
\label{eq:S1S2}g_{0,0}=g_0,\quad \partial_{\la_1}g_{\la_1,0}|_{\la_1=0}=S_1,\quad  \partial_{\la_2}g_{0,\la_2}|_{\la_2=0}=S_2.
\end{equation}
We claim that 
\begin{equation}\label{eq:hess_htop}
\begin{split}
  (\pi_2^* S_1 ,&\mathcal D_{g_0} \tilde \mu(S_2))=
    \int_{SM} (\pi_2^* S_1) (\pi_2^* S_2) \, dm_0 - \frac{2}{n-1}  \partial_{\la_1} |_{\la_1 = 0}\partial_{\la_2}|_{\la_2 =0} h_{\rm top}(g_{\la_1, \la_2})   
    \\&- \int_{SM} \partial_{\la_2}|_{\la_2 = 0} \partial_{\la_1}|_{\la_1 = 0} g_{\la_1, \la_2}(v, v) \, dm_0(v).
    \end{split}
\end{equation}
To see this, we compute $\partial_{\la_2}|_{\la_2=0}\partial_{\la_1}|_{\la_1=0}h_{\mathrm{top}}(g_{\la_1,\la_2})$ using \eqref{eq:KKW}. First, we have
\[
   \partial_{\la_1}|_{\la_1=0} h_{\mathrm{top}}(g_{\la_1,\la_2}) = - \frac{h_{\rm top} (g_{0,\la_2})}{2} \int_{SM} \frac{1}{\Vert v \Vert^2_{{0,\la_2}}} \partial_{\la_1}|_{\la_1=0} g_{\la_1,\la_2} (v, v) \, d \tilde \mu_{0,\la_2}.
    \]
     Second, since $S_2$ has mean zero, using \eqref{eq:KKW} again gives $\partial_{\la_2}|_{\la_2 = 0} h_{\rm top}(g_{0, \la_2}) = 0$. Hence, differentiating the previous equation at $\la_2 = 0$, using $h_{\rm top}(g_0) = n-1$, and rearranging gives \eqref{eq:hess_htop}.

   To prove the desired symmetry and bilinearity, we set $g_{\la_1, \la_2} = g_0 + \la_1 S_1 + \la_2 S_2$ in \eqref{eq:hess_htop}. 
   Then the third term on the right-hand side in \eqref{eq:hess_htop} is zero since the integrand vanishes identically  and the second term is equal, by definition, to $-\tfrac{2}{n-1}d^2_{g_0} h_{\rm top}(S_1, S_2)$. Hence, the left-hand side of \eqref{eq:hess_htop} is bilinear and symmetric, as desired. (A posteriori, this shows the sum of the second and third terms in \eqref{eq:hess_htop} depends only on $S_1$ and $S_2$, in a symmetric and bilinear fashion, for \emph{any} perturbation satisfying \eqref{eq:S1S2}, even though the individual terms may depend on the second variation of $g_{\la_1, \la_2}$.)

   The second step of the proof is to prove the proposition for $S_1 = S_2 = S$.
   To do so, we instead let $g_\la$ be any \emph{volume-preserving} perturbation of $g_0 = g_{\rm hyp}$ with $\dot{g}_0 = S$ (note that $\ddot{g}_0 \neq 0$ here). 
      Then, by \cite[Proposition 5.1.1]{Fla}, the second term on the right-hand side in \eqref{eq:hess_htop}  is given by
$\partial^2_\lambda|_{\lambda=0}h_{\mathrm {top}}(g_\lambda) = {\rm Cov}\left( \pi_2^*T(S), \frac{n-1}{2} \pi_2^*S \right)$, where $T$ is the differential operator in the statement of the proposition. 

To simplify the third term on the right-hand side of \eqref{eq:hess_htop}, we differentiate twice the condition $\mathrm{Vol}_{\la}(M)\equiv \mathrm{Vol}(M)$ (see for instance \cite[Equation $(8.2)$]{Pol}) and use \eqref{eq:pi-trace}. As a result, we obtain
\begin{equation}
\label{eq:Pollicott}
\begin{split}
\int_{SM}\pi_2^*(\ddot g_0)dm_0(v)&=\frac{1}{n\mathrm{Vol}(M)}\int_M \tr_0(\ddot g_0)d\mathrm{vol}_0
\\&=\frac{1}{n\mathrm{Vol}(M)}\|S\|^2_{L^2(M;S^2T^*M)}-\frac{1}{2n\mathrm{Vol}(M)}\| \tr_0(S)\|^2_{L^2(M,\mathrm {vol}_0)}.
\end{split}
\end{equation}

Using \eqref{eq:hess_htop}, \eqref{eq:Pollicott}, and $\partial^2_\lambda|_{\lambda=0}h_{\mathrm {top}}(g_\lambda) = {\rm Cov}\left( \pi_2^*T(S), \frac{n-1}{2} \pi_2^*S \right)$ gives
\begin{align}
\label{eq:quadraticform}
{\rm Cov}&\left( \pi_2^*T(S),  \pi_2^*S \right) = -\Big(-\int_{SM}(\pi_2^*S)^2\,dm_0+\int_{SM}\pi_2^*\ddot g_0\,dm_0+(\pi_2^* S ,\mathcal D_{g_0} \tilde \mu(S))\Big )\\
&=\left(\Vert\pi_2^*S\Vert^2-\frac{1}{n\mathrm{Vol}(M)}\left(\Vert S\Vert^2-\frac{1}{2}\Vert \tr_0(S)\Vert^2\right)-(\pi_2^* S ,\mathcal D_{g_0} \tilde \mu(S))\right),\nonumber
     \end{align}   
which proves the proposition for $S_1 = S_2 = S$.

Finally, we deduce the statement of the proposition by polarization, 
that is, we set $S = S_1 \pm S_2$ in \eqref{eq:quadraticform} and subtract the resulting equations from one another. Since all terms in \eqref{eq:polarization} are symmetric and bilinear, the result follows.
\end{proof}
We now prove the following lemma.
\begin{lemma}
\label{lemm:name}Let $S$ be a symmetric $2$-tensor. Then
    \[{\rm Cov}(\pi_2^*T(\Delta_L S),\pi_2^* S) = {\rm Cov}(\pi_2^*\Delta_L S, \pi_2^*T(S)).\]
\end{lemma}
\begin{remark}
   It can be deduced from the results of \cite[\S 5]{Fla} that $$(S_1,S_2)\mapsto {\rm Cov}(\pi_2^* T(S_1), \pi_2^* S_2)\quad\text{is symmetric.}$$ Since, we will only need to apply it to $S_1=S$ and $S_2=\Delta_LS$, we will not prove this more general statement here.
\end{remark}
\begin{proof}
   We decompose $S = \sum_{i=1}^{+ \infty} S_i$ as in \cite[Lemma 5.3.3]{Fla}. Note that each $S_i$ is an eigentensor of the rough Laplacian, and hence by \eqref{eq:LaplFla} of the Lichnerowicz Laplacian with eigenvalue we will denote by $\mu_i$.  Note that by the definition of $T$ and \eqref{eq:XandD}, we also have for any $i\in \mathbb N$, we have $\pi_2^*T(S_i)$ is cohomologous $\theta_i \pi_2^*S_i$ for some $\theta_i$ depending on $\mu_i.$ Hence, we have ${\rm Cov}(\pi_2^*T(S_i),\cdot)=\theta_i{\rm Cov}(\pi_2^*S_i,\cdot)$. 
   In particular, by \cite[Lemma 5.3.9]{Fla} (see also \cite[Equation (5.3)]{Fla}), we have
   \begin{align*}
   {\rm Cov}&(\pi_2^*T(\Delta_ L S), \pi_2^*S) = \sum_{i,j=1}^{+\infty} {\rm Cov}(\pi_2^*T( \mu_j S_j),\pi_2^* S_i) = \sum_{i,j=1}^{+\infty}\mu_j \theta_j {\rm Cov}( \pi_2^*S_j, \pi_2^*S_i)
   \\&=\sum_{i=1}^{+\infty}\mu_i \theta_i {\rm Cov}( \pi_2^*S_i, \pi_2^*S_i)=\sum_{i,j=1}^{+\infty}\mu_j \theta_i {\rm Cov}(\pi_2^* S_j, \pi_2^*S_i)={\rm Cov}(\pi_2^*\Delta_ L S,\pi_2^*T(S)).
   \end{align*}
   This completes the proof of the lemma.
\end{proof}

\subsection{Lower bound for the Hessian}
In this section, we use our above computations to bound $d^2_{g_0} \Phi(S, S) = \partial^2_{\la =0} \Phi(g_{\la})$ away from zero in terms of $S$ (where $S$ is divergence-free and of mean zero). 
We start by bounding it below by a variance term. 
\begin{proposition}

\label{prop:var}
    Let $(M^n,g_{\rm hyp}) $ be a closed hyperbolic manifold, let $S$ be a divergence-free tensor of zero mean, and let $(g_\la)_{\la \in (-\eps,\eps)}$ be a variation such that $S=\partial_\la|_{\la=0}g_\la$. Then we have
    \begin{align*}
    \partial^2_\la|_{\la=0}\Phi(g_\la)\geq \frac{(n-2)^2}4\mathrm{Var}(\pi_2^*S).
    \end{align*}
\end{proposition}
\begin{proof}

Using Proposition \ref{propHess1}, Proposition \ref{lemma:Flaminio}, and Lemma \ref{lemm:name} we obtain
\[
\partial^2_{\la}|_{\la =0} \Phi(g_{\la}) = 2(n-1) {\rm Cov}(S, T(S)) + {\rm Cov}(\Delta_L S, T(S)).
\]
Now we decompose $S$ into
$S=\sum_{i=1}^{+\infty} S_i, $
where the $S_i$ are eigenvectors of $\nabla^* \nabla$
chosen as in \cite[Lemma 5.3.3]{Fla}.
By \eqref{eq:LaplFla}, the $S_i$ are also eigenvectors of
$\Delta_L$  with eigenvalues we will denote by $\mu_i$. 
Using this, together with the linearity of $T$, the previous equation becomes
\[
\partial^2_{\la}|_{\la =0} \Phi(g_{\la}) = 2(n-1) {\rm Cov}\left(\sum_{i=1}^{+\infty} S_i,\, \sum_{j=1}^{+\infty} T(S_j)\right) + {\rm Cov}\left(\sum_{i=1}^{+\infty} \mu_i S_i,\, \sum_{j=1}^{+\infty} T(S_j)\right).
\]
Using \cite[Lemma 5.3.9]{Fla} yields
\begin{equation*}\label{hess-cov-equality}
\partial^2_\la|_{\la=0}\Phi(g_\la)
=   \sum_{i=1}^{+\infty}(2n-2+\mu_i)\mathrm{Cov}(\pi_2^*S_i,\pi_2^*T(S_i)) .\end{equation*}
Using \eqref{eq:LaplFla}, we see that $\Delta_L\geq 0$ on conformal perturbations. On trace-free perturbations, a Weitzenböck formula shows that $\Delta_L\geq -n$ (see, for instance, \cite[Equation (2.49)]{Hum}). As such, we have $\mathrm{Spec}(\Delta_L|_{\mathrm{Ker}(D_{g_0}^*)})\subset [-n,+\infty)$, which means that we have $2n-2+\mu_i\geq n-2$ for all $i\in\mathbb N$. Thus, using again \cite[Equation (5.3)]{Fla}, we obtain
\begin{align*}
\partial^2_\la|_{\la=0}\Phi(g_\la)
& \geq (n-2)\mathrm{Cov}(\pi_2^*S,\pi_2^*T(S))\geq \frac{(n-2)^2}4\mathrm{Var}(\pi_2^*S), 
\end{align*}
where the last inequality was obtained in \cite[Proposition 1.3.4]{Fla}.
\end{proof}

Next, we use the following \emph{coercive estimate} due to Guillarmou and Lefeuvre on the \emph{generalized X-ray transform} in negative curvature (or equivalently on the variance evaluated on symmetric tensors); see, for instance, \cite[Lemma 2.1]{GuKnLef}. There is $C>0$ such that
\begin{equation*}
%\label{eq:coercive}
\ \forall h\in H^{-1/2}(M;S^2T^*M)\text{ with } \int_{SM}\pi_2^*h \, dm_g=0,\quad \mathrm{Var}_g(\pi_2^*h) \geq C\|\Pi_{\mathrm{Ker}(D_{g}^*)}h\|_{H^{-1/2}}^2,
\end{equation*}
where $\Pi_{\mathrm{Ker}(D_{g}^*)}$ is the orthogonal projection onto solenoidal (or divergence-free) tensors.

Hence, as a direct consequence of Proposition \ref{prop:var}, we obtain the following
\begin{corollary}\label{cor:micro}
      Let $(M^n,g_{\rm hyp}) $ be a closed hyperbolic manifold. Then there exists a constant $C = C(g_0) > 0$ such that, for any divergence-free tensor of zero mean $S$, if we let $(g_\la)_{\la \in (-\eps,\eps)}$ be a variation such that $S=\partial_\la|_{\la=0}g_\la$, we have
    \begin{align*}
    \partial^2_\la|_{\la=0}\Phi(g_\la)\geq C \Vert S \Vert^2_{H^{-1/2}(M)}.
    \end{align*}
\end{corollary}

\section{Taylor expansion of $\Phi$ at a hyperbolic metric}
\label{sec:micro}
In this last section, we deduce Theorem \ref{theo:2} from Corollary \ref{cor:micro}. The argument follows closely that of \cite[Proposition  3.4]{Hum}.
We will consider a second-order Taylor expansion of $\Phi(g)$ near $g = g_{\rm hyp}$. We start with the following control of the remainder term. 

\begin{lemma}\label{lemm:remainder}
    Let $g_{\la} = g_0 + \lambda S$ be a first order perturbation. Then we have
     \[
\Phi(g_{\rm hyp} + S) = \Phi(g_{\rm hyp}) + \partial_{\la}|_{\la=0}\Phi(g_{\la}) + \frac{1}{2} \partial_{\la}^2|_{\la=0}\Phi(g_{\la})+ O( \Vert S \Vert^3_{C^5}).
\]
\end{lemma}

\begin{proof}
The Taylor expansion of the function $F(\la) : = \Phi(g_{\la})$ about $\la = 0$ gives     \[
    |\Phi(g_{\rm hyp} + S) - \Phi(g_{\rm hyp}) - \partial_\la|_{\la=0}\Phi(g_\la) - \frac{1}{2} \partial_{\la}^2|_{\la=0}\Phi(g_\la)| \leq \frac{1}{6}\sup_{\la \in [0,1]} |F'''(\la)|.
    \] 
Now let \[f_{\la}(v) = -\frac{1}{\Vert v \Vert^2_{\la}} {\rm Ric}_{\la}(v, v) + \frac{\mathcal{S}(\la)}{n {\rm Vol}_{\la}(M)},\]
and recall that $F(\la) = (f_{\la}, \tilde \mu_{\la})$, where $( \cdot, \cdot )$ denotes the distributional pairing.
We then have
\begin{equation}\label{eq: 3rd deriv F}
F'''(\la) = (\partial^3_{\la} f_{\la}, \tilde \mu_{\la}) + 3(\partial^2_{\la} f_{\la}, \partial_{\la} \tilde \mu_{\la}) + 3(\partial_{\la} f_{\la}, \partial^2_{\la} \tilde\mu_{\la}) + (f_{\la}, \partial^3_{\la} \tilde \mu_{\la}).
\end{equation}

    Now we use \cite[Theorem C c)]{Con} for the variation of metrics $(g_\la)_{\la \in (-\eps,\eps)}$ and the constant family of potentials identically equal to zero. Then for some small enough $0 < \alpha < 1$, we have that the map
    \begin{align}\label{eq:contreras-map}
    C^{\alpha}(SM) \times \mathcal M_{<0}^{k+1}(M) \to \R, \quad
    (f, g) \mapsto \int_{SM} f d \tilde\mu_g
    \end{align}
    is $C^{k-1}$.
    In particular, the $C^3$-norm of the above map can be controlled by the $C^5$-norm of the variation of metrics. More precisely, for any $j=0,1,2,3$, there exists $C_j>0$ such that
    \begin{equation}
    \label{eq:dual}\forall f\in C^\alpha(SM),\quad |(f,\partial^j_{\la}\tilde \mu_\la)|\leq C_j\|f\|_{C^\alpha}\|S\|_{C^5}^j. 
    \end{equation}
 Indeed, \cite[Theorem C c)]{Con} states that the mapping $g\mapsto \tilde \mu_g \in (C^\alpha(SM))^*$ is a $C^3$ map for a variation of metrics $g_\la$ which is $C^5$. Since our variation is linear, the higher order derivatives of $\tilde \mu_\la$ only depend on $S$ and one has $\|\partial_\la^j\tilde \mu_\la\|_{(C^\alpha(SM))^*}\leq C_j\|S\|_{C^5}^j$ for some constant $C_j>0$ independent of $S$.   The estimate \eqref{eq:dual} then follows from the definition of the dual norm on $(C^\alpha(SM))^*.$ 
 
 In order to conclude the proof of the lemma, we apply \eqref{eq:dual} to $f=\partial_\la^j f_\la$. From the definition of $f_\la$, we see that $g\mapsto f_g$ is a $C^\infty$ function from $C^{2,\alpha}$ metrics to $C^\alpha(SM)$. This means that $\Vert \partial^j_{\la} f_{\la} \Vert_{C^{\alpha}} \leq C\Vert  S \Vert^j_{C^{2, \alpha}} \leq C\Vert  S \Vert_{C^5}^j$, for some constant $C>0$ independent of $S,$ where we used again that the perturbation is linear so the higher derivatives of $f_\la$ only depend on $S$. Combining this last estimate with \eqref{eq:dual}, we obtain $\sup_{\la \in [0,1]} |F'''(\la)| \leq C \Vert S \Vert^3_{C^5}$, which completes the proof.   
\end{proof}

\begin{proof}[Proof of Theorem 1.3]
    Fix $N = N(g_{\rm hyp}) \in \mathbb{N}$ and $\eps = \eps(g_{\rm hyp}) > 0$ to be specified later and let $g$ be any metric such that
    ${\rm Vol}(g) = {\rm Vol}(g_{\rm hyp})$ and 
    $\Vert g - g_{\rm hyp} \Vert_{C^N} < \eps$. 
    By Lemma \ref{slice lemma}, there exists $\phi_g\in \mathrm{Diff}_0(M)$ such that $\phi_g^*g-g_{\rm hyp}\in \mathrm{Ker}(D_{g_{\rm hyp}}^*)$.
  Moreover, there exists $c_g \in \mathbb R$  such that $S: = c_g\phi_g^* g - g_{\rm hyp}$ has zero mean. Note that $S$ is divergence-free as $S = c_g \phi^*g - g_{\rm hyp} = c_g (\phi^* g - g_{\rm hyp}) + (c_g-1) g_{\rm hyp}$ and $g_{\rm hyp}$ is divergence-free. The constant $c_g$ satisfies
$$\int_{M}\tr_{g_{\rm hyp}}(c_g\phi_g^*g-g_{\rm hyp})\,d\mathrm{vol}_{g_{\rm hyp}}=0\iff c_g^{-1}=\frac{1}{n\mathrm{Vol}(M)}\int_M \tr_{g_{\rm hyp}}(\phi_g^*g)\,d\mathrm{vol}_{g_{\rm hyp}}. $$
In particular, the mapping $g\mapsto c_g$ is smooth and, more precisely, by the Lemma \ref{slice lemma}, we have 
\begin{equation}\label{eq:c_g}
    \|g-g_{\mathrm{hyp}}\|_{C^N}\leq \epsilon \ \Rightarrow \ |c_g-1|\leq C\|g-g_{\mathrm{hyp}}\|_{C^N}, 
\end{equation}
for some constant $C>0$ independent of $g.$ Note that if $\eps>0$ is small enough, then $c_g$ is (uniformly in $g$) close to $1.$

Let $g_\la = g_{\rm hyp} + \la S$.  By the definition of $\Phi$ \eqref{eq:Phi2}, we see that for any negatively curved metric $g$ and $c>0$, one has $\Phi(cg)=c^{-1}\Phi(g).$ Then, using Lemma \ref{lemm:remainder}, we have
\[
c_g^{-1} \Phi(g) = \Phi(c_g\phi^*_g g) = \Phi(g_{\rm hyp}) + \partial_\la|_{\la=0}\Phi(g_\la) + \frac{1}{2} \partial_\la^2|_{\la=0}\Phi(g_\la)+ O( \Vert S \Vert^3_{C^5}).
\]
Using that $\Phi(g_{\rm hyp}) = 0$, Lemma \ref{lemm:crit}, and the coercive estimate of Corollary \ref{cor:micro}, we have
\[
c_g^{-1} \Phi(g) \geq C \Vert S \Vert^2_{H^{-1/2}} - C' \Vert S \Vert^3_{C^5},
\]
for some constants $C,C'>0$ independent of the metric $g.$

Since $S$ is small in $C^N$-norm, we can absorb the term $C_n'\|S\|_{C^{5}}^3$ into the left-hand side using a Sobolev embedding and interpolation argument. 
More precisely, for any $\beta>0$, 
$$\|S\|_{C^{5}}^3\leq c\|S\|_{H^{n/2+5+\beta}}^3\leq c'\|S\|_{H^{-1/2}}^2\|S\|_{H^N}\leq c'\|S\|_{H^{-1/2}}^2\|S\|_{C^N}, $$
for any $N>\tfrac 3 2 n+16$ and for some constants $c,c'$ independent of $g$. 

Thus, for any $0<\eps'<\frac{C}{C'c'}$, we have $\Vert S \Vert_{C^N} < \eps'$ implies $\Phi(g) > C \Vert S \Vert^2_{H^{-1/2}}$, for some constant $C$ independent of $g$.
From equation \eqref{eq:c_g}, we deduce that 
$$\|S\|_{C^N}\leq C'\|g - g_0\|_{C^N},$$
for some constant $C'$ independent of $g.$ 
We now choose $\eps$ such that $C' \eps < \eps'$ and this completes the proof.
\end{proof}
\bibliography{references}{}
\bibliographystyle{alpha}
\end{document}